# Recurrences for certain sequences of binomial sums in terms of (generalized) Fibonacci and Lucas polynomials

Johann Cigler

**Abstract.**
We give a simplified presentation of some results about recurrences of certain sequences of binomial sums in terms of (generalized) Fibonacci and Lucas polynomials.

## 1. Introduction

Motivated by the remarkable identities

$$F_n = \sum_{j\in\mathbb{Z}} (-1)^j \binom{n}{\left\lfloor \frac{n+5j+2}{2} \right\rfloor} \tag{1}$$

and

$$F_{n+1} = \sum_{j\in\mathbb{Z}} (-1)^j \binom{n}{\left\lfloor \frac{n+5j}{2} \right\rfloor} \tag{2}$$

for the Fibonacci numbers I obtained in [4] - [9] recurrences for binomial sums of the form

$$A_n(k,m,\ell,z) = \sum_{h\in\mathbb{Z}} \binom{n}{\left\lfloor \frac{n+mh+\ell}{k} \right\rfloor} z^h \tag{3}$$

in terms of (generalized) Fibonacci and Lucas polynomials. Using some results of John P. D'Angelo [2], Eduardo H. M. Brietzke [3] and Dusty E. Grundmeier [10] we give a simplified derivation of these recurrences. As a byproduct we obtain recurrences of subsequences of generalized Fibonacci polynomials. To make the paper accessible to a wider readership, we also provide proofs of some known results that are not part of common knowledge.

## 2. The case k=2

The recurrences of the sequences $\left(A_n(2,m,\ell,z)\right)_{n\geq 0}$ can be expressed using well-known properties of Fibonacci and Lucas polynomials.

The Fibonacci polynomials

$$F_n(x,s) = \sum_{j=0}^{\left\lfloor \frac{n}{2} \right\rfloor} \binom{n-1-j}{j} s^j x^{n-1-2j} \tag{4}$$

satisfy $F_n(x,s) = xF_{n-1}(x,s) + sF_{n-2}(x,s)$ with initial values $F_0(x,s) = 0$ and $F_1(x,s) = 1.$

The Lucas polynomials

$$L_n(x,s) = \sum_{j=0}^{\left\lfloor \frac{n}{2} \right\rfloor} \binom{n-j}{j} \frac{n}{n-j} s^j x^{n-2j} \tag{5}$$

satisfy $L_n(x,s) = xL_{n-1}(x,s) + sL_{n-2}(x,s)$ with initial values $L_0(x,s) = 2$ and $L_1(x,s) = x.$

---

*2020 Mathematics Subject Classification:* 0502, 05A10, 05E05, 11B37, 11B39.
*Key words and phrases:* Fibonacci numbers, Fibonacci and Lucas polynomials, binomial sums, recurrences, elementary symmetric polynomials, Padovan numbers, Perrin numbers, Narayana's cows sequence.

Binet's formulae give

$$F_n(x,s)=\frac{\alpha^n(x,s)-\beta^n(x,s)}{\alpha-\beta} \tag{6}$$

and

$$L_n(x,s)=\alpha^n(x,s)+\beta^n(x,s) \tag{7}$$

with

$$\alpha(x,s)=\frac{x+\sqrt{x^2+4s}}{2},\quad \beta(x,s)=\frac{x-\sqrt{x^2+4s}}{2}. \tag{8}$$

Since

$$\alpha(x+y,-xy)=\frac{x+y+\sqrt{(x+y)^2-4xy}}{2}=x,\quad \beta(x+y,-xy)=\frac{x+y-\sqrt{(x+y)^2-4xy}}{2}=y,$$

we get the well-known formulas

$$L_n(x+y,-xy)=x^n+y^n \tag{9}$$

and

$$F_n(x+y,-xy)=\frac{x^n-y^n}{x-y}. \tag{10}$$

Let $\Delta$ denote the difference operator $\Delta f(x)=f(x+1)-f(x)$ and $E$ the translation operator $Ef(x)=f(x+1)$ on the vector space of polynomials which satisfy $E^i\binom{x}{r}=\binom{x+i}{r}$ and $\Delta^j\binom{x}{r}=\binom{x}{r-j}$.

By (9) we get $L_m(x+1,-x)=x^m+1$ which for $x=\Delta$ gives

$$L_m(E,-\Delta)=\Delta^m+I, \tag{11}$$

where $I$ denotes the identity operator.

Applying (11) to the polynomial $\binom{x}{r}$ gives

$$L_m(E,-\Delta)\binom{x}{r}=\sum_{j=0}^{\left\lfloor \frac{m}{2}\right\rfloor}\binom{m-j}{j}\frac{m}{m-j}(-1)^j\binom{x+m-2j}{r-j}=\binom{x}{r-m}+\binom{x}{r}.$$

Choosing $x=n,\ r=\left\lfloor \frac{n+mh+\ell}{2}\right\rfloor$, multiplying both sides with $z^h$ and summing over all $h\in\mathbb{Z}$ we get

$$\sum_{h\in\mathbb{Z}}\sum_{j=0}^{\left\lfloor \frac{m}{2}\right\rfloor}\binom{m-j}{j}\frac{m}{m-j}(-1)^j\binom{n+m-2j}{\left\lfloor \frac{n-2j+m+m(h-1)+\ell}{2}\right\rfloor}z^h=\sum_{h\in\mathbb{Z}}\binom{n}{\left\lfloor \frac{n+m(h-2)+\ell}{2}\right\rfloor}z^h+\sum_{h\in\mathbb{Z}}\binom{n}{\left\lfloor \frac{n+mh+\ell}{2}\right\rfloor}z^h.$$

Changing the order of summation gives

$$z\sum_{j=0}^{\left\lfloor \frac{m}{2}\right\rfloor}\binom{m-j}{j}\frac{m}{m-j}(-1)^jA_{n+m-2j}(2,m,\ell,z)=\left(1+z^2\right)A_n(2,m,\ell,z).$$

This can be written as

$$\left(zL_m(N,-1)-\left(1+z^2\right)\right)A_n(2,m,\ell,z)=0, \tag{12}$$

if $N$ denotes the operator defined by $N^iA_n=A_{n+i}$.

If $p(N)f(n)=0$ for a polynomial $p(t)$ we call $p(t)$ the *characteristic polynomial of the sequence* $(f(n))$. Note that $p(t)$ is unique up to a multiplicative constant.

**Theorem 1**
*The characteristic polynomial of the sequence* $A_n(2,m,\ell,z)$ *is*

$$p(t)=\left(1+z^2\right)-zL_m(t,-1). \tag{13}$$

.

Since $L_5(x,-1)=x^5-5x^3+5x$ each sequence $A_n(2,5,\ell,-1)$ satisfies the recurrence $\left(N^5-5N^3+5N+2\right)A_n(2,5,\ell,-1)=0.$ Observing that $N^5-5N^3+5N+2=(N+2)\left(N^2-N-1\right)^2$ and $\left(N^2-N-1\right)F_{n+i}=0$ for each $i$ identities (1) and (2) follow from the fact that $A_n(2,5,2,-1)=F_n$ and $A_n(2,5,0,-1)=F_{n+1}$ for $0\le n\le 4$.

**Remark**
The identities (1) and (2) have been found by George Andrews [1]. The proof of the Rogers-Ramanujan identities by Issai Schur [12] contains a $q-$analog of these identities.
For $z=\pm1$ simpler recurrences will be obtained in paragraph 4 with other methods.

The numbers $A_n(2,m,0,-1)$ can be interpreted (cf. e.g. [9]) among other things as the number of the set of all lattice paths in $\mathbb{Z}^2$ which start at the origin, consist of $\left\lfloor \frac{n}{2}\right\rfloor$ north-east steps $U=(1,1)$ and $\left\lfloor \frac{n+1}{2}\right\rfloor$ south-east steps $D=(1,-1)$ and are contained in the strip $-\left\lfloor \frac{m-1}{2}\right\rfloor \le y\le \left\lfloor \frac{m-2}{2}\right\rfloor$. For small $m$ these sequences occur in many contexts ( cf. OEIS[11]):

| $m$ | *OEIS* | first terms |
|---|---|---|
| 3 | *A*000012 | $1,1,1,1,1,1,\cdots$ |
| 4 | *A*016116 | $1,1,2,2,4,4,8,8,\cdots$ |
| 5 | *A*000045 | $1,1,2,3,5,8,13,21,34,55\cdots$ |
| 6 | *A*182522 | $1,1,2,3,6,9,18,27,54,81,162,\cdots$ |
| 7 | *A*028495 | $1,1,2,3,6,10,19,33,61,108,197,\cdots$ |
| 8 | *A*030436 | $1,1,2,3,6,10,20,34,68,116,232,396,\cdots$ |
| 9 | *A*061551 | $1,1,2,3,6,10,20,35,69,124,241,440,846,\cdots$ |

For $z=1$ and $\ell=0$ we get

| $m$ | *OEIS* | first terms |
|---|---|---|
| 2 | *A*000079 | $1,2,4,8,16,32,\cdots$ |
| 3 | *A*001045 | $1,1,3,5,11,21,43,85,\cdots$ |
| 4 | *A*011782 | $1,1,2,4,8,16,32,64,128,\cdots$ |
| 5 | *A*099163 | $1,1,2,3,7,12,27,49,106,199,\cdots$ |
| 6 | *A*005578 | $1,1,2,3,6,11,22,43,86,171,342,683,\cdots$ |

$A_{2n}(2,2m,0,1)=\sum_{h\in\mathbb{Z}}\binom{2n}{n+hm}$ can be interpreted as the number of closed walks of length $2n$ on a vertex of the cyclic graph on $2m$ nodes or equivalently as the number of lattice paths with steps $U$ and $D$ from $(0,0)$ to $(2n,2hm)$ for some $h\in\mathbb{Z}.$

$A_{2n+1}(2,2m,0,1)=\sum_{h\in\mathbb{Z}}\binom{2n+1}{n+hm}$ is the number of walks of length $2n+1$ between two adjacent vertices.

Therefore, we get $A_{2n+2}(2,2m,0,1)=2A_{2n+1}(2,2m,0,1).$

## 3. The general case

### 3.1. Generalized Fibonacci and Lucas polynomials.

The Fibonacci polynomials $F_n^{(k)}(x,s)$ are defined by $F_n^{(k)}(x,s)=xF_{n-1}^{(k)}(x,s)+sF_{n-k}^{(k)}(x,s)$ with initial values $F_0^{(k)}(x,s)=0$ and $F_n^{(k)}(x,s)=x^{n-1}$ for $0<n<k.$

The Lucas polynomials $L_n^{(k)}(x,s)$ are defined by $L_n^{(k)}(x,s)=xL_{n-1}^{(k)}(x,s)+sL_{n-k}^{(k)}(x,s)$ with initial values $L_0^{(k)}(x,s)=k$ and $L_n^{(k)}(x,s)=x^n$ for $0<n<k.$

From $\left(1-xz-sz^k\right)\sum_{n\ge0}F_n^{(k)}(x,s)z^n=\left(1-xz-sz^k\right)\left(z+xz^2+\cdots+x^{k-1}z^k+\cdots\right)=z$ we get

$$\sum_{n\ge0}F_{n+1}^{(k)}(x,s)z^n=\frac{1}{1-xz-sz^k} \tag{14}$$

$$\frac{1}{1-xz-sz^k}=\sum_{\ell\ge0}\left(xz+sz^k\right)^\ell=\sum_{j,\ell}\binom{\ell}{j}s^jz^{jk}x^{\ell-j}z^{\ell-j}=\sum_{j(k-1)+\ell=n}z^n\sum_{j=0}^{\ell}\binom{\ell}{j}s^jx^{\ell-j}=\sum_{n\ge0}\sum_{j=0}^{\left\lfloor\frac{n}{k}\right\rfloor}\binom{n-(k-1)j}{j}s^jx^{n-kj}$$

gives

$$F_{n+1}^{(k)}(x,s)=\sum_{j=0}^{\left\lfloor\frac{n}{k}\right\rfloor}\binom{n-(k-1)j}{j}s^jx^{n-kj}. \tag{15}$$

From

$$\left(1-xz-sz^k\right)\sum_{n\ge0}L_n^{(k)}(x,s)z^n=\left(1-xz-sz^k\right)\left(k+xz+x^2z^2+\cdots+x^{k-1}z^{k-1}+x^kz^k+\cdots\right)=k-(k-1)xz$$

we get

$$\sum_{n\ge0}L_n^{(k)}(x,s)z^n=\frac{k-(k-1)xz}{1-xz-sz^k}. \tag{16}$$

Thus we get $L_n^{(k)}(x,s)=kF_{n+1}^{(k)}(x,s)-(k-1)xF_n^{(k)}(x,s).$

Since

$$k\binom{n-(k-1)j}{j}-(k-1)\binom{n-1-(k-1)j}{j}=\binom{n-(k-1)j}{j}+(k-1)\binom{n-1-(k-1)j}{j-1}$$

$$=\binom{n-(k-1)j}{j}\left(1+(k-1)\frac{j}{n-(k-1)j}\right)=\binom{n-(k-1)j}{j}\frac{n}{n-(k-1)j}$$

we get for $n\ge k$

$$L_n^{(k)}(x,s)=\sum_{j=0}^{\lfloor \frac{n}{k}\rfloor}\binom{n-(k-1)j}{j}\frac{n}{n-(k-1)j}s^j x^{n-kj}. \tag{17}$$

Other generalizations of $F_n(x,s)$ are the polynomials $G_n^{(k)}(x,s)$ with generating function

$$\sum_{n\geq 0} G_n^{(k)}(x,s)z^n=\frac{1}{1-s^{k-2}xz^{k-1}-s^{k-1}z^k}. \tag{18}$$

They satisfy $G_n^{(k)}(x,s)=s^{k-2}xG_{n-1}^{(k)}(x,s)+s^{k-1}G_{n-k}^{(k)}(x,s)$ with initial values $G_0^{(k)}(x,s)=1,$ $G_n^{(k)}(x,s)=0$ for $0<n\leq k-2$ and $G_{k-1}^{(k)}(x,s)=s^{k-2}x$ and

$$G_n^{(k)}(x,s)=\sum_{j=\left\lceil\frac{(k-2)n}{k}\right\rceil}^{\frac{(k-1)n}{k}}\binom{n-j}{(k-1)n-kj}s^j x^{(k-1)n-kj}. \tag{19}$$

**Proof**

$$\frac{1}{1-s^{k-2}xz^{k-1}-s^{k-1}z^k}=\sum_{\ell\geq 0}\left(s^{k-2}z^{k-1}\right)^{\ell}\left(x+sz\right)^{\ell}=\sum_{\ell\geq 0, i\geq 0}\binom{\ell}{i}s^{(k-2)\ell+i}x^{\ell-i}z^{(k-1)\ell+i}$$

If we set $(k-2)\ell+i=j$ and $(k-1)\ell+i=n$ then $\ell=n-j$ and $(k-2)(n-j)+i=j,$ which gives

$n-j-i=n-j-\left(j-(k-2)(n-j)\right)=n-2j+(k-2)n-(k-2)j=(k-1)n-kj.$

This gives

$$\sum_{\ell\geq 0, i\geq 0}\binom{\ell}{\ell-i}s^{(k-2)\ell+i}x^{\ell-i}z^{(k-1)\ell+i}=\sum_{n\geq 0}z^n\sum_j\binom{n-j}{(k-1)n-kj}s^j x^{(k-1)n-kj},$$

which is (19).

Since $\binom{n}{k}=0$ for $k<0$ we can also write

$$G_n^{(k)}(x,s)=\sum_{j=0}^{n}\binom{n-j}{(k-1)n-kj}s^j x^{(k-1)n-kj}=\sum_{j=0}^{n}\binom{n-j}{(k-1)j-(k-2)n}s^j x^{(k-1)n-kj}$$
$$=\sum_{j=0}^{n}\binom{n-j}{(k-2)(j-n)+j}s^j x^{(k-1)n-kj}=\sum_{j=0}^{n}\binom{j}{n-(k-1)j}s^{n-j}x^{kj-n}. \tag{20}$$

The analog of the Lucas polynomials are the polynomials

$$H_n^{(k)}(x,s)=\sum_{j=\left\lceil\frac{(k-2)n}{k}\right\rceil}^{\frac{(k-1)n}{k}}\binom{n-j}{(k-1)n-kj}\frac{n}{n-j}s^j x^{(k-1)n-kj}. \tag{21}$$

They can also be written as

$$H_n^{(k)}(x,s)=\sum_{j=0}^{n-1}\binom{n-j}{(k-1)n-kj}\frac{n}{n-j}s^j x^{(k-1)n-kj}=\sum_{j=1}^{n}\binom{j}{n-(k-1)j}\frac{n}{j}s^{n-j}x^{kj-n}.$$

The initial values are $H_0^{(k)}(x,s)=k,$ $H_n^{(k)}(x,s)=0$ for $0<n<k-1,$ $H_{k-1}^{(k)}(x,s)=(k-1)s^{k-2}x.$

$$(1-s^{k-2}xz^{k-1}-s^{k-1}z^k)\sum_{n\geq 0}H_n^{(k)}(x,s)z^n=(1-s^{k-2}xz^{k-1}-s^{k-1}z^k)\left(k+(k-1)s^{k-2}xz^{k-1}+\cdots\right)=k-s^{k-2}xz^{k-1}$$

implies

$$\sum_{n\geq 0}H_n^{(k)}(x,s)z^n=\frac{k-s^{k-2}xz^{k-1}}{1-s^{k-2}xz^{k-1}-s^{k-1}z^k}. \tag{22}$$

**Remark**

The numbers $G_n^{(2)}(1,1)$ are the Fibonacci numbers $F_{n+1}$. By (18) the numbers $\left(G_n^{(3)}(1,1)\right)_{n\geq 0} = \left(1,0,1,1,1,2,2,3,4,5,7,9,12,\cdots\right)$ are the Padovan numbers OEIS [11], A000931.
For $4 \leq k \leq 7$ the numbers are listed in OEIS, A017817, A017827, A017837, A017847.
The numbers $H_n^{(2)}(1,1)$ are the Lucas numbers $L_n$. By (22) the numbers $\left(H_n^{(3)}(1,1)\right)_{n\geq 0} = \left(3,0,2,3,2,5,5,7,10,12,17,22,29,39,51,\cdots\right)$ are the Perrin numbers A001608.
The numbers for $4 \leq k \leq 7$ occur in A050443, A087937, A087936, A306755.

### 3.2 Recurrences for the general case

The method for $k=2$ can be generalized to give

**Lemma 2**

*Let*

$$\sum_j c(m,k,j)x^j(1+x)^{i_j m-kj} = 0, \tag{23}$$

*where* $i_j$ *is an integer such that* $i_j m - kj \geq 0$. *Then*

$$\sum_j c(m,k,j)z^{i_j} A_{n+i_j m-kj}\left(k,m,\ell,z\right) = 0. \tag{24}$$

**Proof**

Applying the operator $\sum_j c(m,k,j)\Delta^j E^{i_j m-kj}$ to $\binom{x}{\left\lfloor \frac{n+mh+\ell}{k} \right\rfloor}$ we get

$$\sum_j c(m,k,j)\Delta^j E^{i_j m-kj}\binom{x}{\left\lfloor \frac{n+mh+\ell}{k} \right\rfloor} = \sum_j c(m,k,j)\binom{x+i_j m-kj}{\left\lfloor \frac{n+mh+\ell-kj}{k} \right\rfloor} = 0.$$

If we set $x=n$ this gives

$$\sum_j c(m,k,j)\binom{n+i_j m-kj}{\left\lfloor \frac{n+i_j m-kj+m(h-i_j)+\ell}{k} \right\rfloor} = 0.$$

Multiplying with $z^h$ and summing over $h \in \mathbb{Z}$ gives

$$\sum_j c(m,k,j)\sum_{h\in\mathbb{Z}} z^h \binom{n+i(j)m-kj}{\left\lfloor \frac{n+i_j m-kj+m(h-i_j)+\ell}{k} \right\rfloor} = \sum_j c(m,k,j)z^{i_j} A_{n+i_j m-kj}\left(k,m,\ell,z\right) = 0.$$

Let us first consider the case $k=1$.

Since $-1+\sum_{j=0}^{m}\binom{m}{j}(-1)^j x^j (1+x)^{m-j} = 0$ is of the form (23) we get

$$z\sum_{j=0}^{m}(-1)^j\binom{m}{j}A_{n+m-j}\left(1,m,\ell,z\right) = A_n\left(1,m,\ell,z\right). \tag{25}$$

For $k>1$ the following identity of the form (23) has been found in [4]:
There exist uniquely determined integers $a_{m,k,j}$ such that

$$(1+x)^m-1+\sum_{i=1}^{k-1}\sum_{j=\left\lceil\frac{(i-1)m}{k-1}\right\rceil}^{\left\lfloor\frac{im}{k}\right\rfloor}a_{m,k,j}x^j(x+1)^{im-kj}=0. \tag{26}$$

Note that 1, $x^m$, and $(1+x)^m$ are of the form $x^j(x+1)^{im-kj}$.

Suppose first that such a formula exists. Then the polynomial $\sum_{j=1}^{m}a_{m,k,j}z^j$ has the root $\omega_m^{-k}(\omega_m-1)$, where $\omega_m$ denotes a primitive $m-$th root of unity. The most obvious polynomial with this root is

$$\prod_{i=1}^{m}\left(z-\omega_m^{-kj}\left(\omega_m^j-1\right)\right)=\sum_{j=1}^{m}b_{m,k,j}z^j. \tag{27}$$

This led in [4], Lemma 7.1 to the formula

$$(-1)^{k(m-1)}\sum_{j=1}^{m}b_{m,k,j}x^j(1+x)^{\operatorname{mod}(-kj,m)}=(x+1)^m-1, \tag{28}$$

where $\operatorname{mod}(\ell,m)$ denotes the least non-negative residue of $\ell$ modulo $m$.

For example for $(k,m)=(3,4)$ we get $\sum_{j=1}^{4}b_{4,3,j}z^j=-4z-2z^2+z^4$. This gives
$4x(1+x)+2x^2(1+x)^2-x^4=4x(1+x)^{4-3}+2x^2(1+x)^{8-2\cdot3}-x^4(1+x)^{12-4\cdot3}=(x+1)^4-1.$

From another point of view John P. D'Angelo [2] showed that there is a uniquely determined polynomial

$$f_{m,k}(x,y)=1-\prod_{j=0}^{m-1}\left(1-\omega_m^jx-\omega_m^{kj}y\right) \tag{29}$$

which satisfies the following 4 conditions:

1) $f(0,0)=0,$
2) $f(x,y)=1$ when $x+y=1,$
3) $\deg f=m,$
4) $f\left(\omega_m x,\omega_m^k y\right)=f(x,y),$

It turns out that $r(m,k,x,y)=x^m+(-1)^{k(m-1)}\sum_{j=1}^{m}b_{m,k,j}(-y)^jx^{\operatorname{mod}(-kj,m)}$ satisfies these conditions: 1),3),4) are obvious and 2) follows from (28).

This observation led to a simpler approach of (26) in [8], which I will now present in a slightly different form which emphasizes the analogy with the case $k=2$.
Let $\alpha_1(x,s),\alpha_2(x,s),\cdots,\alpha_k(x,s)$ be the roots of the polynomial $a_{k,1}(z,x,s)=z^k-xz^{k-1}-s$ and let

$$a_{k,m}(z,x,s)=\prod_{j=1}^{k}\left(z-\alpha_j^m(x,s)\right)=\sum_{i=0}^{k}(-1)^ie_{i,k,m}(x,s)z^{k-i}. \tag{30}$$

To compute the elementary symmetric polynomials $e_{i,k,m}(x,s)$ observe that

$$\prod_{i=1}^{k}\left(1-\alpha_i^mz^m\right)=\prod_{i=1}^{k}\prod_{j=0}^{m-1}\left(1-\omega_m^j\alpha_iz\right)=\prod_{j=0}^{m-1}\prod_{i=1}^{k}\left(1-\omega_m^j\alpha_iz\right)=\prod_{j=0}^{m-1}\left(1-\omega_m^jxz-\omega_m^{kj}sz^k\right)=\Phi_{k,m}(z,x,s).$$

Therefore

$$\prod_{j=0}^{m-1}\left(1-\omega_m^j xz-\omega_m^{kj} sz^k\right)=\prod_{i=1}^{k}\left(1-\alpha_i^m z^m\right)=\sum_{i=0}^{k}(-1)^i e_{i,k,m}(x,s) z^{im}. \tag{31}$$

Here $e_{0,k,m}(x,s)=1$ and

$$e_{k,k,m}(x,s)=(-1)^{(k-1)m}s^m, \tag{32}$$

because $(-1)^k\alpha_1\alpha_2\cdots\alpha_k=-s$ and therefore $e_{k,k,m}(x,s)=(\alpha_1\alpha_2\cdots\alpha_k)^m=(-1)^{(k-1)m}s^m$. For $m=0$ we get from (30)

$$e_{i,k,0}=\binom{k}{i}. \tag{33}$$

The left-hand side of (31) is a linear combination of $s^j z^{kj} x^{m-j} z^{m-j}=s^j x^{m-j} z^{m+(k-1)j}$. Only such terms can occur where $m+(k-1)j=im$ for some $i$. The polynomials $e_{i.k.m}(x,s)$ are linear combinations of $s^j x^{m-j}=s^j x^{im-kj}$. We have $0\le j\le m=in-kj+j$ and $im-kj\ge 0$. This gives $j\le\frac{im}{k}$ and $j=kj-(i-1)m\ge 0$, i.e. $j\ge\frac{(i-1)m}{k}$.

Therefore for $1\le i<k$

$$p_i(m,k,x,s)=(-1)^{i+1}e_{i,k,m}(x,s)=\sum_{j=\left\lceil\frac{(i-1)m}{k}\right\rceil}^{\left\lfloor\frac{im}{k}\right\rfloor}a(m,k,j)s^j x^{im-kj} \tag{34}$$

for some integers $a(m,k,j)$. The left-hand side of (31) vanishes for $(x+s,z)=(1,1)$, because the factor for $j=0$ vanishes. This implies

$$\Phi_{k,m}(1,x+1,-x)=\sum_{i=0}^{k}(-1)^i e_{i,k,m}(x+1,-x)=0, \tag{35}$$

or equivalently

$$\sum_{i=1}^{k-1}p_i(m,k,x+1,-x)=1+(-1)^{k(m-1)}x^m. \tag{36}$$

By Lemma 2 we get

$$\sum_{i=1}^{k-1}z^i p_i(m,k,N,-1)A_n(k,m,\ell,z)=\left(1+(-1)^{k(m-1)}z^k\right)A_n(k,m,\ell,z). \tag{37}$$

**Theorem 3**

*The characteristic polynomial of the sequence* $\left(A_n(k,m,\ell,z)\right)_{n\ge 0}$ *is*

$$p_{m,k,z}(t)=\left(1+(-1)^{k(m-1)}z^k\right)-\sum_{i=1}^{k-1}z^i p_i(m,k,t,-1). \tag{38}$$

As shown by Dusty Grundmeier [10] the characteristic polynomials $P_{i,k}(t,x,s)$ of the sequences $\left(p_i(m,k,x,s)\right)_{m\ge 0}$ have degree $\binom{k}{i}$.

For $i=1$ and $i=k-1$ they are explicitly given by

$$\begin{aligned}P_{1,k}(t,x,s)&=t^k-xt^{k-1}-s,\\ P_{k-1,k}(t,x,s)&=t^k+(-1)^{k-1}s^{k-2}xt-s^{k-1}.\end{aligned} \tag{39}$$

In order to verify (39) let us recall Newton's identities
$p_m - p_{m-1}e_1 + p_{m-2}e_2 - \cdots + (-1)^{m-1} p_1 e_{m-1} + (-1)^m m s_m$
for the power sums $p_m(x_1,\cdots,x_k) = \sum_i x_i^m$ in terms of the elementary symmetric polynomials $e_i(x_1,\cdots,x_k)$.

Since $p_1(m,k,x,s) = e_{1,k,m}(x,s) = \sum_{i=1}^{k} \alpha_i(x,s)^m$ and $e_{i,k,1}(x,s)$ is the $i-$th elementary symmetric polynomial of $\alpha_1,\cdots,\alpha_k$, we get $e_{1,k,m}(x,s) = xe_{1,k,m-1}(x,s) + se_{1,k,m-k}(x,s)$.
From $\prod_{j=0}^{m-1}\left(1-\omega_m^j xz\right) = 1 - x^m z^m$ we get from (31) that for $0 < m < k$ the initial values of $e_{1,k,m}(x,s)$ are $x^m$.
This gives for $m > 0$

$$p_1(m,k,x,s) = L_m^{(k)}(x,s). \tag{40}$$

Note that $e_{k-1,k,m}(x,s) = \left(\alpha_1\alpha_2\cdots\alpha_k\right)^m \left(\frac{1}{\alpha_1^m} + \cdots + \frac{1}{\alpha_k^m}\right)$ and let

$$u(z) = \prod_{i=1}^{k}\left(\alpha_1\alpha_2\cdots\alpha_k z - \alpha_i\right) = \prod_{i=1}^{k}\left((-1)^{k-1}sz - \alpha_i\right) = s^k z^k - (-1)^{k-1} x s^{k-1} z^{k-1} - s.$$

Then $-\frac{1}{s} z^k u\left(\frac{1}{z}\right) = z^k + (-1)^{k-1} s^{k-2} xz - s^{k-1}$ has the roots $\frac{\alpha_1\alpha_2\cdots\alpha_k}{\alpha_i}$.
Therefore $e_{k-1,k,m}(x,s)$ has the characteristic polynomial $p(t) = t^k + (-1)^{k-1} s^{k-2} xt - s^{k-1}$
This gives with Newton's identities (cf. [4] and [10])

$$\sum_{m\geq 0} p_{k-1}(m,k,x,s) z^m = \frac{k - (-s)^{k-2} x z^{k-1}}{1 - (-s)^{k-2} x z^{k-1} - s^{k-1} z^k} \tag{41}$$

Therefore, we have

$$p_{k-1}(m,k,x,s) = H_m^{(k)}\left(x,(-1)^k s\right). \tag{42}$$

In the general case no explicit formulas for $p_i(m,k,x,s)$ with $2 \leq i < k-1$ are known.
Consider for example $p_2(m,4,x,s)$. Its characteristic polynomial is $t^6 + st^4 + sx^2t^3 - s^2t^2 - s^3$.
The first terms of $\left(p_2(m,4,x,s)\right)_{m\geq 0}$ are
$(-6,0,2s,3sx^2,-6s^2,-5s^2x^2,2s^3-3s^2x^4,14s^3x^2,-6s^4+8s^3x^4,-18s^4x^2+3s^3x^6,2s^5-25s^4x^4,\cdots)$.
By (34) we get for $i \in \{0,1\}$

$$\begin{aligned} p_2\left(4m+2i,4,x,s\right) &= \sum_{j=0}^{m+i} a\left(4m+2i,4,2m+i-j\right) s^{2m+i-j} x^{4j}, \\ p_2\left(4m+1+2i,4,x,s\right) &= \sum_{j=0}^{m+i} a(4m+1+2i,4,2m+i-j) s^{2m+i-j} x^{4j+2}. \end{aligned} \tag{43}$$

Computations suggest the following values for the first coefficients:

$$a(4m,4,2m)=-6,\ \ a(4m,4,2m-1)=\frac{16m}{m-2}\binom{m+1}{4},$$
$$a(4m,4,2m-2)=-\frac{16m\left(4m^2-15\right)}{(m+3)(m-3)(m-4)}\binom{m+3}{8},$$
$$a(4m,4,2m-3)=\frac{32m\left(8m^2-35\right)}{(m-5)(m-6)(m-7)}\binom{m+4}{12},$$
$$a(4m+2,4,2m+1)=2,\ \ a(4m+2,4,2m+1-1)=-\frac{(4m+2)^2}{(m-1)(m-2)}\binom{m+1}{4},$$
$$a(4m+2,4,2m+1-2)=\frac{4\left(4m+2\right)^2}{(m-3)(m-4)}\binom{m+3}{8},$$
$$a(4m+2,4,2m+1-3)=\frac{2\left(4m+2\right)^2\left(-105+8m+8m^2\right)}{\left(m-4\right)\left(m-5\right)\left(m-6\right)\left(m-7\right)}\binom{m+4}{12},$$

$$a(4m+1,4,2m)=-m(4m+1),\ \ a(4m+1,4,2m-1)=\frac{2(4m+1)(4m-3)}{(m-2)(m-3)}\binom{m+2}{6},$$
$$a(4m+1,4,2m-2)=-\frac{8(4m+1)\left(4m^2+9m-10\right)}{(m-4)(m-5)(m-6)}\binom{m+3}{10}$$

$$a(4m+3,4,2m+1)=(m+1)(4m+3),\ \ a(4m+3,4,2m+1-1)=-\frac{2(4m+3)(4m+7)}{(m-2)(m-3)}\binom{m+2}{6},$$
$$a(4m+3,4,2m+1-2)=\frac{8(m+4)(4m+3)\left(4m^2-m-15\right)}{(m-3)(m-4)(m-5)(m-6)}\binom{m+3}{10}.$$

Perhaps someone can find a general formula.

**Some examples**

For $k=3, z=-1$ the characteristic polynomial is

$$p_{m,3,-1}(t)=\left(1+(-1)^m\right)+p_1(m,3,t,-1)-p_2(m,3,t,-1). \qquad (44)$$

Here we get $\left(A_n(3,1,0,-1)\right)_{n\ge0}=\left(0,0,0,\cdots\right)$ with $p_{1,3,-1}(t)=t$,
$\left(A_n(3,2,0,-1)\right)_{n\ge0}=\left(1,-2,2,0,-4,8,-8,0,16,-32,32,0,-64,128,-128,0,256,-512,\cdots\right)$ with $p_{2,3,-1}(t)=2+2t+t^2$,
$\left(A_n(3,3,0,-1)\right)_{n\ge0}=\left(0,0,0,\cdots\right)$ with $p_{1,3,-1}(t)=t^3$,
$\left(A_n(3,4,0,-1)\right)_{n\ge0}=\left(0,0,1,-2,-2,8,-6,-20,48,0,-164,232,232,-1120,792,2704,-6528,0,\cdots\right)$ with $p_{1,4,-1}(t)=2-4t+2t^2+t^4$,
$\left(A_n(3,5,0,-1)\right)_{n\ge0}=\left(1,0,0,1,-5,0,5,-30,25,25,-175,275,0,-1000,2250,-1375,-5000,\cdots\right)$ with $p_{1,5,-1}(t)=5t-5t^2+t^5$.

## 4. Simpler recurrences for $k=2$ and $z=\pm1$.

**4.1.** By (13) the sequence $\left(A_n(2,m,\ell,-1)\right)_{n\ge0}$ has characteristic polynomial $L_m(t,-1)+2$.
The polynomials $L_n(x,-1)+2$ have non-trivial factors.

$$L_{2n}(x,-1)+2=L_n(x,-1)^2, \qquad (45)$$

$$L_{2n+1}(x,-1)+2=(x+2)\left(F_{n+1}(x,-1)-F_n(x,-1)\right)^2. \tag{46}$$

To prove these identities observe that $\alpha(x,-1)\beta(x,-1)=1$, $(\alpha(x,-1)-\beta(x,-1))^2=x^2-4$, $(\alpha(x,-1)-1)^2=(x-2)\alpha(x,-1)$, $(\beta(x,-1)-1)^2=(x-2)\beta(x,-1)$, $(\alpha(x,-1)-1)(\beta(x,-1)-1)=2-x$.

(45) follows from

$$L_m(x,-1)^2=\left(\alpha(x,-1)^m+\beta(x,-1)^m\right)^2=\alpha(x,-1)^{2m}+2\left(\alpha(x,-1)\beta(x,-1)\right)^m+\beta(x,-1)^{2m}=L_{2m}(x,-1)+2.$$

To verify (46) observe that

$$\left(\alpha-\beta\right)^2\left(F_{m+1}(x,-1)-F_m(x,-1)\right)^2=\left(\alpha^{m+1}-\beta^{m+1}-\alpha^m+\beta^m\right)^2=\left(\alpha^m(\alpha-1)-\beta^m(\beta-1)\right)^2$$
$$=\alpha^{2m}(\alpha-1)^2+\beta^{2m}(\beta-1)^2-2\alpha^m\beta^m(\alpha-1)(\beta-1)=\alpha^{2m}(x-2)\alpha+\beta^{2m}(x-2)\beta-2(2-x)$$
$$=(x-2)\left(\alpha^{2m+1}+\beta^{2m+1}+2\right).$$

**4.2.** Identities (45) and (46) suggest

**Theorem 4**

$$L_m(N,-1)A_n(2,2m,\ell,-1)=0 \tag{47}$$

*and*

$$\left(F_{m+1}(N,-1)-F_m(N,-1)\right)A_n(2,2m+1,\ell,-1)=0 \tag{48}$$

*for each* $\ell\in\mathbb{Z}$.

In these cases Lemma 2 is not applicable because the operators depend on $m$ and not on $2m$ or $2m+1$. Let us first consider (47). Here the method of Lemma 2 works.

Applying the identity $L_m(E,-\Delta)=\Delta^m+I$ to $\binom{x}{r}$ gives

$$\sum_{j=0}^{\left\lfloor \frac{m}{2}\right\rfloor}(-1)^j\binom{m-j}{j}\frac{m}{m-j}\binom{x+m-2j}{r-j}=\binom{x}{r}+\binom{x}{r-m}.$$

Setting $x=n$ and $r=\left\lfloor\frac{n+\ell-2mh+m}{2}\right\rfloor$, multiplying the identity by $(-1)^h$ and summing over $h\in\mathbb{Z}$ we get

$$\begin{aligned}&\sum_{j=0}^{\left\lfloor \frac{m}{2}\right\rfloor}(-1)^j\binom{m-j}{j}\frac{m}{m-j}\sum_{h\in\mathbb{Z}}(-1)^h\binom{n+m-2j}{\left\lfloor\frac{n+\ell-2mh-2j+m}{2}\right\rfloor}\\&=\sum_{h\in\mathbb{Z}}(-1)^h\binom{n}{\left\lfloor\frac{n+\ell-2mh+m}{2}\right\rfloor}+\sum_{h\in\mathbb{Z}}(-1)^h\binom{n}{\left\lfloor\frac{n+\ell-2m(h+1)+m}{2}\right\rfloor}\end{aligned} \tag{49}$$

Observing that $\sum_{h\in\mathbb{Z}}(-1)^h\binom{n+m-2j}{\left\lfloor\frac{n+m-2j+\ell-2mh}{2}\right\rfloor}=A_{n+m-2j}\left(2,2m,\ell,-1\right)$

and that the right-hand side of (49) vanishes we get (47).

For the proof of (48) we recall the proof of [5] which uses an idea of E. Brietzke [3]. Consider the numbers

$$t(n,k)=(-1)^k\binom{n}{\left\lfloor \frac{n+k}{2}\right\rfloor}. \tag{50}$$

They satisfy $t(n,k)=-t(n-1,k-1)-t(n-1,k+1)$ with $t(0,0)=1, t(0,1)=-1$ and $t(0,k)=0$ for all other $k\in\mathbb{Z}$.

Let $s(n,k)$ on $\mathbb{N}\times\mathbb{Z}$ be the function which satisfies the same recurrence, but with initial values $s[0,k)=[k=0]$. Writing in row $n$ the values $f(n,k)$, we get a signed Pascal's triangle:

| | | | | | | | | |
|---|---|---|---|---|---|---|---|---|
| 0 | 0 | 0 | 0 | 1 | 0 | 0 | 0 | 0 |
| 0 | 0 | 0 | −1 | 0 | −1 | 0 | 0 | 0 |
| 0 | 0 | 1 | 0 | 2 | 0 | 1 | 0 | 0 |
| 0 | −1 | 0 | −3 | 0 | −3 | 0 | −1 | 0 |
| 1 | 0 | 4 | 0 | 6 | 0 | 4 | 0 | 1 |

Then the function $t(n,k)$ looks like

| | | | | | | | | |
|---|---|---|---|---|---|---|---|---|
| 0 | 0 | 0 | 0 | 1 | −1 | 0 | 0 | 0 |
| 0 | 0 | 0 | −1 | 1 | −1 | 1 | 0 | 0 |
| 0 | 0 | 1 | −1 | 2 | −2 | 1 | −1 | 0 |
| 0 | −1 | 1 | −3 | 3 | −3 | 3 | −1 | 1 |
| 1 | −1 | 4 | −4 | 6 | −6 | 4 | −4 | 1 |

Define linear operators $N$ and $K$ for functions on $\mathbb{N}\times\mathbb{Z}$ by

$$\begin{aligned} Nf(n,k)&=f(n+1,k),\\ Kf(n,k)&=f(n,k-1). \end{aligned} \tag{51}$$

They give

$$s(n,k)=Ns(n-1,k)=-s(n-1,k-1)-s(n-1,k+1)=-\left(K+K^{-1}\right)s(n-1,k)=(-1)^n\left(K+K^{-1}\right)^n s(0,k).$$

The initial values give $t(0,k)=(1-K)s(0,k)$ and linearity gives $t(n,k)=(1-K)s(n,k).$

Let $\mathcal{F}$ be the linear space of all functions on $\mathbb{N}\times\mathbb{Z}$ which are finite linear combinations of functions $K^j s(n,k)$ for $j\in\mathbb{Z}.$

For $f\in\mathcal{F}$ we have

$$Nf=-\left(K+K^{-1}\right)f. \tag{52}$$

For $x\to -x,\ y\to -\frac{1}{x}$ in (10) we get $F_n\left(-x-\frac{1}{x},-1\right)=(-1)^{n-1}\sum_{j=0}^{n-1}x^{n-1-2j}$ and

$F_{n+1}\left(-x-\frac{1}{x},-1\right)-F_n\left(-x-\frac{1}{x},-1\right)=(-1)^n\sum_{j=-n}^{n}x^j.$

Therefore on $\mathcal{F}$ we have

$$\left(F_{m+1}(N,-1)-F_m(N,-1)\right)=(-1)^m\sum_{j=-m}^{m}K^j. \tag{53}$$

This implies $\left(F_{m+1}(N,-1)-F_m(N,-1)\right)s(0,k)=(-1)^m[|k|\le m]$ and

$$\left(F_{m+1}(N,-1)-F_m(N,-1)\right)\sum_{h\in\mathbb{Z}} s(n,\ell-(2m+1)h)=(-1)^m.$$

Since $t(n,k)=(1-K)s(n,k)$ we also have

$\left(F_{m+1}(N,-1)-F_m(N,-1)\right)\sum_{h\in\mathbb{Z}} t(n,\ell-(2m+1)h)=0.$ (48) follows from

$$\sum_{h\in\mathbb{Z}} t(n,\ell-(2m+1)h)=(-1)^\ell\sum_{h\in\mathbb{Z}}(-1)^h\binom{n}{\left\lfloor \frac{n+\ell-(2m+1)h}{2}\right\rfloor}=(-1)^\ell A_n(2,2m+1,\ell,-1).$$

**4.3.** For $z=1$ there are also simpler recurrences.

**Theorem 5**

$$(N-2)F_m(N,-1)A_n\left(2,2m,\ell,1\right)=0 \tag{54}$$

*and*

$$\left(L_{m+1}(N,-1)-L_m(N,-1)\right)A_n\left(2,2m+1,\ell,1\right)=0. \tag{55}$$

**Remark**

Note that

$$L_{2m}(x,-1)-2=(x^2-4)F_m(x,-1)^2 \tag{56}$$

and

$$L_{2m+1}(x,-1)-2=\frac{\left(L_{m+1}(x,-1)-L_m(x,-1)\right)^2}{x-2}. \tag{57}$$

Formula (56) follows from

$\left(\alpha(x,-1)^m-\beta(x,-1)^m\right)^2=\alpha(x,-1)^{2m}+\beta(x,-1)^{2m}-2\left(\alpha(x,-1)\beta(x,-1)\right)^m.$

Formula (57) follows from

$$\left(L_{m+1}(x,-1)-L_m(x,-1)\right)^2=\left(\alpha^{m+1}+\beta^{m+1}-\alpha^m-\beta^m\right)^2=\left(\alpha^m(\alpha-1)+\beta^m(\beta-1)\right)^2$$

$$=\alpha^{2m}(\alpha-1)^2+\beta^{2m}(\beta-1)^2+2\alpha^m\beta^m(\alpha-1)(\beta-1)=\alpha^{2m}(x-2)\alpha+\beta^{2m}(x-2)\beta+2(2-x)$$

$$=(x-2)\left(\alpha^{2m+1}+\beta^{2m+1}-2\right)$$

**Proof of Theorem 5.**

The function $u(n,\ell)=\binom{n}{\left\lfloor \frac{n+\ell}{2}\right\rfloor}$ satisfies $u(n,\ell)=u(n-1,\ell-1)+u(n-1,\ell+1)$ and

$u(0,\ell)=\left[\ell\in\{0,1\}\right].$ Therefore $u(n,\ell)=(1+K)S(0,\ell),$ where $S(n,\ell)$ satisfies the same recurrence as $u(n,\ell)$ but with initial values $S(0,\ell)=[\ell=0].$

Let $\mathcal{G}$ be the linear space of all finite linear combinations $K^jS(n,\ell)$ for $j\in\mathbb{Z}.$ On $\mathcal{G}$ we have $N=K+K^{-1}.$

The identity

$$\left(x+\frac{1}{x}-2\right)F_m\left(x+\frac{1}{x},-1\right)(1+x)=\frac{(x-1)^2}{x}\left(\frac{x^m-\frac{1}{x^m}}{x-\frac{1}{x}}\right)(1+x)=(x-1)\left(x^m-\frac{1}{x^m}\right)=\frac{1}{x^m}-\frac{1}{x^{m-1}}-x^m+x^{m+1}$$

implies

$(E-2)F_m(E,-1)u(0,\ell)=\left(K^m-K^{m-1}-K^{-m}+K^{-m-1}\right)S(0,\ell)$

and therefore

$$(E-2)F_m(E,-1)\sum_{h\in\mathbb{Z}} K^{2hm}u(0,\ell)=\sum_{h\in\mathbb{Z}} K^{2hm}\left(K^m-K^{m-1}-K^{-m}+K^{-m-1}\right)S(0,\ell)=0.$$

This gives (54).
From

$$\left(L_m\left(x+\frac{1}{x},-1\right)-L_{m-1}\left(x+\frac{1}{x},-1\right)\right)(1+x)=(1+x)\left(x^m+\frac{1}{x^m}-x^{m-1}-\frac{1}{x^{m-1}}\right)=\frac{1}{x^m}-\frac{1}{x^{m-2}}-x^{m-1}+x^{m+1}$$

we get (55).

## 5. Some final remarks

It is perhaps interesting to note that (30) also gives the recurrences of subsequences of the generalized Fibonacci polynomials $F_n^{(k)}(x,s)$.
By (30) the characteristic polynomials of the subsequences $\left(F_{mn+r}^{(k)}(x,s)\right)$ are given by

$$q_{m,k}(t)=\sum_{i=0}^{k}(-1)^i e_{i,k,m}(x,s)t^{k-i}. \tag{58}$$

For $k=3$ this reduces to

$$q_{m,3}(t)=t^3-e_{1,3,m}(x,s)t^2+e_{2,3,m}(x,s)t-e_{3,3,m}(x,s). \tag{59}$$

By (40)

$$e_{1,3,m}(x,s)=\sum_{j=0}^{\left\lfloor\frac{m}{3}\right\rfloor}\binom{m-2j}{j}\frac{m}{m-2j}s^j x^{m-3j}. \tag{60}$$

The first terms are $3,\ x,\ x^2,\ x^3+3s,\ x^4+4sx,\ x^5+5sx^2,\ x^6+6sx^3+3s^2,\cdots$.
By (42)

$$e_{2,3,m}(x,s)=\sum_{j=\left\lceil\frac{m}{3}\right\rceil}^{\left\lfloor\frac{2m}{3}\right\rfloor}(-1)^j\binom{m-j}{2m-3j}\frac{m}{m-j}s^j x^{2m-3j}. \tag{61}$$

The first terms are $3,\ 0,\ -2sx,\ 3s^2,\ 2s^2x^2,\ -5s^3x,\ 3s^4-2s^3x^3,\ 7s^4x^2,\cdots$.
By (32)

$$e_{3,3,m}=s^m. \tag{62}$$

**Remark**
The sequence $\left(G_n\right)_{n\geq 0}=\left(1,1,1,2,3,4,6,9,13,19,\cdots\right)$ with $G_n=F_{n+1}^{(3)}(1,1)$ is called Narayana's cows sequence, OEIS [11], A000930.
From (59) we get the recurrences $G_{m(n+3)}=a_mG_{m(n+2)}-b_mG_{m(n+1)}+G_{mn}$, where $a_m=e_{1,3,m}(1,1)$ satisfies $a_m=a_{m-1}+a_{m-3}$ with initial values $a_0=3,a_1=1,a_2=1$ and $b_m=e_{2,3,m}(1,1)$ satisfies $b_m=-b_{m-2}+b_{m-3}$ with initial values $b_0=3,b_1=0,b_2=-2$.
The first terms are $\left(a_m\right)_{m\geq 1}=\left(1,1,4,5,6,10,15,21,\cdots\right)$ and $\left(b_m\right)_{m\geq 1}=\left(0,-2,3,2,-5,1,7,-6,\cdots\right)$.
For example, we get $G_{2n}=G_{2n-2}+2G_{2n-4}+G_{2n-6}$ and $G_{3n}=4G_{3n-3}-3G_{3n-6}+G_{3n-9}$.

Fakultät für Mathematik, Universität Wien
johann.cigler@univie.ac.at